\documentclass[11pt]{article}
\usepackage{amsmath}
\usepackage{amssymb}
\def\k{\mathbf{k}}
\newcommand{\Mod}[1]{\ (\mathrm{mod}\ #1)}
\newtheorem{theorem}{\hspace*{\parindent}Theorem}
\newtheorem{lemma}{\hspace*{\parindent}Lemma}
\newtheorem{corollary}{\hspace*{\parindent}Corollary}
\newtheorem{definition}{\hspace*{\parindent}Definition}
\title{Hypergeometric representations and differential-difference relations for some kernels appearing in mathematical physics}

\author{Dmitrii B.\:Karp$^{\rm a}$\footnote{Corresponding author.  E-mail: dmitriibkarp@tdtu.edu.vn}, Yuri B. Melnikov$^{\rm b}$~~and Irina V. Turuntaeva$^{\rm c}$
\\[10pt]\small{\textit{$\phantom{1}^a$Ton Duc Thang University, Ho Chi Minh City, Vietnam}}
\\
\small{\textit{$\phantom{1}^b$Agentur Kronstadt GmbH (Hamburg, Germany) and Oxford Progress Ltd (Oxford, UK)}}
\\
\small{\textit{$\phantom{1}^c$Far Eastern Federal University, Vladivostok/Ussurijsk, Russia}}
}

\date{}

\begin{document}

\maketitle

\bigskip
\begin{center}
\parbox{12cm}{
\small\textbf{Abstract.} The paper is an investigation of the analytic properties of a new class of special functions that appear in the kernels of a class of integral operators underlying the dynamics of matter relaxation processes in attractive fields.  These functions, recently introduced by the second author, generate the kernels of the principal parts of these operators and play an important role in understanding their spectral characteristics.  We reveal the representations of  these functions in terms of the Gauss and Clausen hypergeometric functions and present differential-difference and differential equations they satisfy.  Mathematically, the results include calculation of certain trigonometric double integrals and derivation of their other properties. Furthermore, they represent a potentially useful tool in matter relaxation in an external field, the study of nanoelectronic electrolyte-based systems and dynamics of charge carriers in media with obstacles.}
\end{center}

\bigskip

Keywords: \emph{generalized hypergeometric function, integral operator, differential-difference equation, integral representation, matter relaxation in a field}

\bigskip

MSC2010: 47G10, 45A05, 33C20, 33C05

\section{Introduction and background}

In the present paper we continue the discussion of the special functions $\Xi_N^{[\bf k]}(x)$ introduced in \cite{MelnikovJMP2016} and partially investigated in \cite{MelnikovJMP2016,MelnikovJCAM2018} (mostly for the case $N = 1$). These functions have naturally appeared as kernels of a particular class of integral operators $K_\varphi$
upon their generalisation to higher dimensional phase spaces \cite{MelnikovJMP2016}. These operators are defined on the Hilbert space $L_2(\Omega)$, $\Omega\subset{\bf R}^M$, as follows
$$
K_\varphi: u({\bf x})\,\mapsto\,\int_\Omega\frac{u({\bf x})\varphi({\bf s})-u({\bf s})\varphi({\bf x})}{|{\bf x}-{\bf s}|}\,d^M{\bf s}\,\,.
$$
Usually, $\Omega$ is a bounded convex domain with a smooth boundary. For collision operators the symbol $\varphi$ is known as the {\it equilibrium distribution function} \cite{MelnikovJMP2016,MelnikovLMP1998,Melnikov1999,MelnikovJMP2001,Melnikov2004} and is interpreted as the probability density. For the matter relaxation processes \cite{MelnikovJMP2016,MelnikovJMP2017-1,MelnikovJMP2017-2} the symbol $\varphi$ plays the role of {\it an external field}.  Typically, the class of operators $K_\varphi$ is restricted by the requirements that $\varphi({\bf x})$ is an {\it acceptable function} \cite{MelnikovJMP2016}, i.e. positive, smooth (i.e.  Lipschitz-1), summable (belongs to $L_1(\Omega)$) and uniformly separated from zero in $\Omega$. The last requirement first appeared for mathematical reasons (as its violation in dimension $1$ impedes obtaining the important spectral estimates  \cite{MelnikovJMP2016,Melnikov2004}). However, it also has a physical background: if we consider a conductor with an induced field, this field, in general, vanishes nowhere in the body of the conductor (a remark due to P.\:Avtonomov).
\par
The operators  $K_\varphi$ were initially introduced as collision operators in some models of non-equilibrium statistical physics  \cite{AntoniouTasaki,PetroskyOrdonez,Prigogine}. Their rigorous mathematical study in dimension 1 started in late 1990-s \cite{MelnikovLMP1998,Melnikov1999,MelnikovJMP2001,Melnikov2004,MelnikovYarevsky}, and has been continued recently in \cite{Loesche,MelnikovJMP2016}. Spectral analysis of these operators has demonstrated encouraging links with various fields of mathematics \cite{Melnikov2004,MelnikovJMP2016}. Another important application in physics has been traced  in \cite{MelnikovJMP2016,MelnikovJMP2017-2}: the operators $K_\varphi$ determine the dynamics of a matter relaxation process in an external attractive field $\varphi({\bf x})$.  Such processes in nanostructures  are of particular interest to physicists due to the need to model electrolyte relaxation in thin films \cite{KYB2006}. Similar spectral problems occur also in the models of the electronic transport through resistive multichannel quantum wires \cite{Gantmakher,Knobchen,Kunze1994,Kunze1995-1,Kunze1995-2,Lenk1994,Lenk1998}.
\par
The paper is organised as follows. In Section 2 we refer to the previous results and introduce notation and definitions. In Sections 3 and 4 we investigate the functions $\Xi_1^{[k]}(x)$ and $\Xi_2^{\bf{[k]}}(x)$, respectively and prove the key Theorem~\ref{th:Xi2representation}. Section 5 is devoted to conclusions and further tasks.

\section{Preliminaries and definitions}

A recent generalisation of the operator $K_\varphi$ to higher dimensions $M \ge 2$  \cite{MelnikovJMP2016} led to the introduction of a new class of special functions \cite[(64)]{MelnikovJMP2016}:
\begin{equation}\label{eq:Xi-defined}
\Xi_N^{[{\bf k}]}(x):=\frac{1}{2N}\int_0^{\pi}\cdots\int_0^{\pi}\,\frac{\prod_{m=1}^N\cos{k_m \theta_m}}{\bigl[1+x^2-2xN^{-1}\sum_{l=1}^N \cos{\theta_l}\bigr]^{1/2}}
d\theta_1\cdots{d\theta_N}.
\end{equation}
Here $0\le{x}<1$ and the multi-index ${\bf k}\in{\bf Z}^N_+$. We call $N=1,2,\ldots$ {\it the rank} and ${\bf k} \in {\bf Z}^N_+$ {\it the (multi-) order} of the function $\Xi_N^{[{\bf k}]}(x)$.  They were introduced and briefly considered in \cite{MelnikovJMP2016}, and studied in a greater detail in \cite{MelnikovJCAM2018} in dimension $2$ ($M=N+1=2$).
These functions (for $N=1$ and $N=2$) form the main object of the present investigation.
\par
The operator $K_\varphi$ is the Friedrichs extension \cite{AkhiezerGlazman1966,BirmanSolomyak1987} of the core operator defined on ${\cal D} = C^1(\Omega)\cap{L_1(\Omega)}\subset L_2(\Omega)$, where $\Omega$ is a bounded domain with smooth boundary. As shown in  \cite{Melnikov2004,MelnikovJMP2016}, it is self-adjoint in the weighted space ${\cal H}_\varphi = L_2\bigl(\Omega,\,d^N{\bf x}/\varphi({\bf x})\bigr)$.
\par
Let us note, that in $1$-dimensional case $M=1$ the operator $K_\varphi$ is not an integral operator, as the factor $|x - s|^{-1}$ in the integrand generates a non-integrable singularity in a vicinity of the point $s=x$. On the other hand, for any $M \ge 2$ this singularity is integrable and the  operator $K_\varphi$ can be decomposed  as \cite{MelnikovJMP2016}
$$
K_\varphi: \,u({\bf x})\,\mapsto\,u({\bf x})\,\int_\Omega\frac{\varphi({\bf s})}{|{\bf x}-{\bf s}|}\,d^M{\bf s}\,
-\,\varphi({\bf x})\,\int_\Omega\frac{u({\bf s})}{|{\bf x}-{\bf s}|}\,d^M{\bf s},
$$
where $d^M{\bf s}$ denotes the $M$-dimensional Lebesgue measure. Hence, for $M \ge 2$ it is the sum of the integral operator with the kernel
$$
\kappa({\bf x},{\bf s}):=\,-\,\varphi({\bf x})\,|{\bf x}-{\bf s}|^{-1}
$$
and the operator of multiplication by the function
$$
q_\varphi({\bf x})\,:=\,\int_\Omega\frac{\varphi({\bf s})}{|{\bf x}-{\bf s}|}\,d^M{\bf s}.
$$
\par
For $M\ge2$, the known results in the spectral analysis of the generalized reference operator $K_0$ are incomplete \cite{MelnikovJMP2016} (contrary to the one-dimensional case), which prevents the application of the direct methods used in \cite{Melnikov2004} for the analysis of the operators $K_\varphi$.   Moreover, the calculation of the quadratic form $\kappa_\varphi[u] := \langle K_\varphi u, u\rangle_{L_2(\Omega; d{\bf x}/\varphi({\bf x}))}$ is not obvious for $M\ge2$. $K_0$ here stands for the operator $K_\varphi$ with the symbol $\varphi({\bf x})\equiv1$ in $\Omega$, which is inconsistent with the notation $K_{\varphi}$, but is inherited from the nomenclature used in the early paper \cite{MelnikovLMP1998} and the subsequent  works.
\par
Still, some results in this direction can be obtained under certain specific conditions discussed below. See also \cite{MelnikovJMP2016}. Namely, let us define the so-called {\it completely admissible domain}:
\begin{definition}\label{df:pseudotorus}
A bounded domain $\Omega \subset \mathbb{R}^M$, $M=N+1\ge 2$ is called completely admissible (or pseudo-torus) if there exists a change of coordinates such that each point of $\Omega$ can be written as ${\bf x} = \{r, {\bf\Theta}\}$, where $r\in[0,R)$ and ${\bf\Theta}:=(\theta_1,\theta_2,\ldots,\theta_N)$, $\theta_k\in[0,2\pi)$, are called the pseudo-toroidal coordinates. Hence, $\Omega$ can be represented  as
\begin{equation}\label{eq1:Subsec4.1}
\Omega=\widehat{\bf T}^N:=[0,R)\times\underbrace{\mathbb{S}^1\times\mathbb{S}^1\times\cdots\times\mathbb{S}^1}_{N\,{\rm times}},
\end{equation}
where $\mathbb{S}^1$ denotes the circle.
\end{definition}
Note that in dimension $M=2$ (i.e. for $N=1$) the simplest example of the pseudo-torus is the annulus $\widehat{\bf T}^1=\{{\bf x}:L_0\le|{\bf x}|\le{L}\}$, $r=R(L-L_0)^{-1}(|{\bf x}|-L_0)$. In the degenerate case $L_0=0$ it reduces to the disk  $\widehat{\bf T}^1_{\rm deg} = \{{\bf x}:0\le|{\bf x}|\le{R}\}$, $r=|{\bf x}|$. In dimension $M = 3$ (i.e. for $N=2$) the simplest pseudo-torus is the spherical ring $\widehat{\bf T}^2$ with the outer radius $L$, the inner radius $L_0$ and $R=L-L_0$. In the degenerate case $L_0=0$ the pseudo-torus $\widehat{\bf T}^2_{\rm deg}$ is topologically equivalent to a three-dimensional ball. These cases are considered in Sections 3 and 4 below.
\par
As shown in \cite{MelnikovJMP2016}, the operator $K_\varphi$ defined by the symbol $\varphi({\bf x})$ having a specific symmetry in the $N+1$-dimensional pseudo-torus $\Omega_R^{[N+1]}$ allows for the reduction of the corresponding spectral problem for the operator $K_\varphi$ on ${\Omega}$ to a set of independent one-dimensional spectral problems for the partial operators $K_\varphi^{[{\bf k}]}$, ${\bf k}\in{\bf Z}^{N}_+$. In the present paper we consider the cases $N = 1$ and $N = 2$.  For $N = 1$ we take $\Omega_R^{[2]}$ to be the disk $\Omega_R^{[2]}:=\{{\bf x}:|{\bf x}| \le R\}$. Any point ${\bf x}\in\Omega_R^{[2]}\subset{\bf R}^2$ can be written in polar coordinates  as ${\bf x}=(r,\theta)$, $0\le r:=|{\bf x}|\le R$, $0\le\theta<2\pi$. The requirement to be imposed on the function $\varphi({\bf x})$ is its rotation invariance: $\varphi({\bf x})=\varphi(r)$. For $N=2$ the domain is the degenerate torus $\Omega_R^{[3]}:=\{{\bf x}: \mathrm{dist}({\bf x},\mathcal{C})<R\}$, where $\mathcal{C}$ is the (one-dimensional) circle of radius $R$ centered at the origin. Any point ${\bf x}\in\Omega_R^{[3]}\subset{\bf R}^3$ can be written as ${\bf x}=(r, {\bf \Theta})$, where $r$ is the distance from the circle
 $\mathcal{C}$ and ${\bf \Theta} = (\theta_1, \theta_2)$, $0\le\theta_1,\theta_2<2\pi$. The requirement on the function $\varphi({\bf x})$ is independence of $\theta_1$, $\theta_2$: $\varphi({\bf x})=\varphi(r)$. It was shown in \cite{MelnikovJMP2016} that under the above assumptions the Hilbert spaces
$$
{\cal H}^{[{\bf k}]} := \bigl\{u({\bf x})\in L_2({\bf R}^{N}): u({\bf x})=e^{i\langle{\bf k},{\bf\Theta}\rangle}\,u^{[{\bf k}]}(r)\bigr\},
$$
where $e^{i\langle{\bf k},{\bf\Theta}\rangle}:= \prod_{l=1}^N \, e^{i k_l \theta_l}$, are invariant subspaces of the operator $K_\varphi$, i.e. $K_\varphi: {\cal H}^{[{\bf k}]}\to{\cal H}^{[{\bf k}]}$. This means that $K_\varphi$ can be decomposed into the orthogonal sum  of the partial operators
$$
K_\varphi=\sum_{{\bf k}\in{\bf Z}_+^N}\oplus K_\varphi^{[{\bf k}]}
$$
acting in ${\cal H}^{[{\bf k}]}$. It was also demonstrated in \cite{MelnikovJMP2016} that these partial operators $K_\varphi^{[{\bf k}]}$ act as follows:
\begin{equation}\label{eq:partialK}
K_\varphi^{[{\bf k}]}:  e^{i\langle{\bf k},{\bf\Theta}\rangle}\,u^{[{\bf k}]}(r)\mapsto e^{i\langle{\bf k},{\bf\Theta}\rangle}\,\frac{2N-1}{2^N\,N!\,\sqrt{N}}\,\int_0^R \rho^{2N-1}\,\bigl[u(r)\varphi(\rho)Z_N^{[{\bf 0}]}(r,\rho)-u(\rho)\varphi(r)Z_N^{[{\bf k}]}(r,\rho)\bigr]\,d\rho,
\end{equation}
where the kernel $Z_N^{[{\bf k}]}(r, \rho)$ can be expressed in terms of the functions $\Xi_N^{[{\bf k}]}(x)$  defined in \eqref{eq:Xi-defined}.

\section{Rank-1 functions $\Xi_1^{[k]}(x)$}

In our previous papers \cite{MelnikovJMP2016,MelnikovJCAM2018} we used a simplified notation for the case $N=1$, i.e. $\xi_k(x):=\Xi_1^{[k]}(x)$.
In the present paper, where both cases $N=1$ and $N=2$ will be considered, we keep the general notation $\Xi_N^{[{\bf k}]}$ to avoid confusion.

According to \cite[(18.12.11)]{NIST} we have
\begin{equation}\label{eq:LegendreGenerating}
\frac{1}{\sqrt{1+x^2-2xz}}=\sum_{n=0}^{\infty}P_n(z)x^n,
\end{equation}
where $P_n$ denotes the Legendre polynomial \cite[p.305]{AAR}, \cite[section~18.3]{NIST}.  Substituting this expansion into (\ref{eq:Xi-defined})  leads to
$$
\Xi_1^{[k]}(x) = \sum_{n=0}^{\infty}\tau^{[k]}_{n}x^n
$$
with
$$
\tau^{[k]}_{n}=\frac{1}{\pi}\int_{0}^{\pi}P_n(\cos(\theta))\cos(k\theta)d\theta.
$$
On the other hand, according to $\lambda=1/2$ case of \cite[(18.5.11)]{NIST},
$$
P_n(\cos(\theta))=\sum_{j=0}^{n}\frac{(1/2)_j(1/2)_{n-j}}{j!(n-j)!}\cos((n-2j)\theta),
$$
where $(a)_k=\Gamma(a+k)/\Gamma(a)$ is the rising factorial.  Hence, we have
$$
\tau^{[k]}_{n}=\frac{1}{\pi}\int_{0}^{\pi}P_n(\cos(\theta))\cos(k\theta)d\theta
=\frac{1}{\pi}\sum_{j=0}^{n}\frac{(1/2)_j(1/2)_{n-j}}{j!(n-j)!}\int_{0}^{\pi}\cos((n-2j)\theta)\cos(k\theta)d\theta.
$$
In view of the orthogonality relation \cite[(4.26.10-11)]{NIST}
\begin{equation}\label{eq:cos-ortho}
\int_{0}^{\pi}\cos((n-2j)\theta)\cos(k\theta)d\theta=\left\{\!\!\!\begin{array}{l}0,~~~n-2j\ne{k},\\\pi/2,~~n-2j=k\ne0,\\\pi,~~n-2j=k=0,\end{array}\right.
\end{equation}
this amounts to
$$
\tau^{[k]}_{n}=\frac{1}{2}\frac{(1/2)_{(n-k)/2}(1/2)_{(n+k)/2}}{((n-k)/)!((n+k)/2)!},~~\text{where}~n\ge{k}~\text{and}~n-k~\text{is even}.
$$
If $n-k$ is odd or $n<k$ then $\tau^{[k]}_{n}=0$.  Finally, changing $n-k$ to $2m$ we get
\begin{multline*}
\Xi_1^{[k]}(x)=\frac{1}{2}\sum_{\stackrel{n=k}{n-k~\text{is even}}}^{\infty}\frac{(1/2)_{(n-k)/2}(1/2)_{(n+k)/2}}{((n-k)/2)!((n+k)/2)!}x^n
=\frac{1}{2}\sum_{m=0}^{\infty}\frac{(1/2)_{m}(1/2)_{m+k}}{m!(m+k)!}x^{2m+k}
\\
=\frac{(1/2)_kx^k}{2k!}\sum_{m=0}^{\infty}\frac{(1/2)_{m}(1/2+k)_{m}}{(1+k)_mm!}x^{2m}.
\end{multline*}
In the second equality we have used the easily verifiable identity $(a)_{m+k}=(a)_k(a+k)_m$.
Thus we have demonstrated our first theorem.

\begin{theorem}\label{th:Xi1via2F1}
The rank-1 function $\Xi_1^{[k]}$ can be expressed as
\begin{equation}\label{eq:RankOneFinal}
\Xi_1^{[k]}(x)=\frac{(1/2)_kx^k}{2k!}{_{2}F_1}\!\left(\!\!\begin{array}{l}1/2,1/2+k\\k+1\end{array}\vline\,\,x^2\!\right),
\end{equation}
where ${_{2}F_1}$ is the Gauss hypergeometric functions whose main properties can be found in \cite[Chapter~15]{NIST}.
\end{theorem}
Representation (\ref{eq:RankOneFinal}) permits the derivation of the analytic and asymptotic (in the neighborhood of $x=1$, in particular)  properties of the functions $\Xi_1^{[k]}(x)$ from the multitude of the known facts about the Gauss hypergeometric function ${_{2}F_1}$.  In particular, we get

\begin{corollary}\label{cr:rank1diffeq}
The function $\Xi^{[k]}_1(x)$ satisfies the second-order differential equation
\begin{equation}\label{eq:Xi1diffeq}
\biggl[x^2(1-x^2)\frac{d^2}{dx^2}+x(1-3x)\frac{d}{dx}+x^2(k^2-1)-k^2\biggr]\Xi^{[k]}_1(x)=0.
\end{equation}
\end{corollary}
\textbf{Proof.}  The hypergeometric function ${_{2}F_1}$ from (\ref{eq:RankOneFinal}) satisfies the hypergeometric differential equation (see \cite[(2.3.5)]{AAR}, \cite[2.1(1)]{HTF1}, \cite[15.10.1]{NIST},  or \cite[(2.1.6)]{SlavyanovLay2000})
$$
\Bigl[z(1-z)D_z^2+\bigl(k+1-(k+2)z\bigr)D_z-\bigl(1/4+k/2\bigr)\Bigr]{_{2}F_1}\left(1/2,1/2+k; k+1;z\right)=0.
$$
We take $z=x^2$ here and apply the identities
$$
D_z=\frac{d}{dz}=\frac{1}{2x}\frac{d}{dx},~~~~D_z^2=\frac{1}{4x^2}\frac{d^2}{dx^2}-\frac{1}{4x^3}\frac{d}{dx}
$$
to the above differential equation to obtain ($D_x=d/dx$):
$$
\Bigl[\frac{1-x^2}{4}D_x^2+\frac{1}{2x}(k+1/2-(k+3/2)x^2)D_x-2k-1\Bigr]{_{2}F_1}\left(1/2,1/2+k; k+1;x^2\right)=0.
$$
Then, employing the formulas
\begin{align*}
&D_x[g(x)]=x^{-\alpha}\left\{D_x-\alpha{x^{-1}}\right\}x^{\alpha}g(x),
\\[5pt]
&D_x^2[g(x)]=x^{-\alpha}\left\{D_x^2-2\alpha{x^{-1}}D_x+\alpha(\alpha+1)x^{-2}\right\}x^{\alpha}g(x)
\end{align*}
with $\alpha=k$ and $g(x)={_{2}F_1}\left(1/2,1/2+k; k+1;x^2\right)$ after tedious but elementary calculations we arrive at (\ref{eq:Xi1diffeq}). $\hfill\square$

Following Fuchs' terminology \cite[Appendix~F.3]{AAR}, \cite[section~1.1]{SlavyanovLay2000}, the point $x$ of equation (\ref{eq:Xi1diffeq}) is called {\it singular}
if the functional coefficient at $d^2/dx^2$ vanishes at $x$, and (in our case) the point $x=\infty$.  Hence, there are three such points: $x=0$, $x=1$ and $x=\infty$.
According to \cite[Appendix~F.3]{AAR}, \cite[Chapter~2]{SlavyanovLay2000}, all these singular points are regular.

We also recall the differential-difference relation discovered in \cite[Lemma~1]{MelnikovJCAM2018} which can also be re-derived from (\ref{eq:RankOneFinal}).
\begin{corollary}\label{cr:rank1dif-dif}
The function $\Xi^{[k]}_1(x)$ satisfies the following differential-difference relation\emph{:}
$$
\biggl[x\frac{d}{dx}+\frac{1}{2}\biggr]\left(\Xi_1^{k}(x)+\Xi_1^{k+2}(x)\right)=\biggl[(1+x^2)\frac{d}{dx}+x\biggr]\Xi_1^{k+1}(x).
$$
\end{corollary}

\section{Rank-2 functions $\Xi_2^{\bf{[k]}}(x)$}

For $N=2$ functions $\Xi_2^{{\bf [k]}}$ are defined by formula (\ref{eq:Xi-defined}) with $N=2$, where ${\bf k}=(k_1,k_2)$.  Note that permutation symmetry
$\Xi_2^{(k_1,k_2)}=\Xi_2^{(k_2,k_1)}$ allows us to restrict our attention to the case $0\le{k_1}\le{k_2}$ without loss of  generality.  The main result of this section is the following theorem.

\begin{theorem}\label{th:Xi2representation}
 Suppose $0\le{k_1}\le{k_2}$ are integers, $\k=(k_1,k_2)$.  Put $s=\lfloor(k_2-k_1)/2\rfloor$.  Then
\begin{multline}\label{eq:theta2-3F2}
\Xi^{[{\bf k]}}_2(x)=\frac{\pi^2(1/2)_{k_1+k_2}(2x)^{k_1+k_2}}{4^{k_1+k_2+1}k_1!k_2!}\sum_{j=0}^{s}\frac{((k_1-k_2)/2)_j((k_1-k_2+1)/2)_j(k_1+k_2+1/2)_{2j}x^{2j}}{(k_1+1)_j(k_2+1)_{j}(k_1+k_2+1)_{j}j!}
\\
\times{_{3}F_2}\!\left(\!\!\begin{array}{l}1/2+j,k_1+1/2+j,k_1+k_2+1/2+2j\\k_2+1+j,k_1+k_2+1+j\end{array}\vline\,\,-x^2\!\right),
\end{multline}
where ${_{3}F_2}$ denotes Clausen's generalized hypergeometric function \emph{\cite[Chapter~16]{NIST}}.
Note that $(0)_0=1$ in the above formula.
\end{theorem}

If $s=0$ in the above theorem the function $\Xi^{[{\bf k]}}_2$ simplifies as follows.
\begin{corollary}\label{cr:Xi2smallk}
Suppose in addition to  the assumptions of Theorem~\ref{th:Xi2representation} that $k_2-k_1\le 1$. Then
\begin{equation}\label{eq:corollaryXi2}
\Xi^{[{\bf k]}}_2(x)=\frac{\pi^2(1/2)_{k_1+k_2}(2x)^{k_1+k_2}}{4^{k_1+k_2+1}k_1!k_2!}
{_{3}F_2}\!\left(\!\!\begin{array}{l}1/2,k_1+1/2,k_1+k_2+1/2\\k_2+1,k_1+k_2+1\end{array}\vline\,\,-x^2\!\right).
\end{equation}
\end{corollary}

The next lemma will play a crucial role in the proof of Theorem~\ref{th:Xi2representation}.
\begin{lemma}\label{lm:Al}
For given non-negative integers $l$, $k_1$, $k_2$ define
\begin{equation}\label{eq:A-defined}
A_l(k_1,k_2):=\iint\limits_{0\le\theta_1,\theta_2\le\pi}\!\!\!(\cos(\theta_1)+\cos(\theta_2))^{l}\cos(k_1\theta_1)\cos(k_2\theta_2)d\theta_1d\theta_2.
\end{equation}
Then
\begin{equation}\label{eq:A-computed}
A_l(k_1,k_2)\!=\!\left\{\!\!\!\begin{array}{lr}
0, & l\not\equiv{k_1+k_2}\Mod{2}~\text{\emph{or}}~l<k_1+k_2
\\[10pt]
\displaystyle\binom{k_1+k_2}{k_1}\dfrac{(\pi^2/2^l)(k_1+k_2+1)_{2N}^2}{(k_1+1)_N(k_2+1)_N(k_1+k_2+1)_NN!}, &l\equiv{k_1+k_2}\Mod2~\text{\emph{and}}~l\ge{k_1+k_2},
\end{array}\right.
\end{equation}
where $N=(l-k_1-k_2)/2$ and
$$
\binom{p}{s}:=\frac{p(p-1)(p-2)\cdots(p-s+1)}{s!}
$$
is the binomial coefficient.
\end{lemma}

\textbf{Proof.} First, we will need the formula
\begin{multline*}
[\cos(\theta)]^p=\left(\frac{e^{i\theta}+e^{-i\theta}}{2}\right)^p=\frac{1}{2^p}\sum_{s=0}^{p}\binom{p}{s}e^{i\theta{s}}e^{-i\theta(p-s)}
=\frac{1}{2^p}\sum_{s=0}^{p}\binom{p}{s}[\cos(\theta(2s-p))+i\sin(\theta(2s-p))]
\\
=\frac{1}{2^p}\sum_{s=0}^{p}\binom{p}{s}\cos(\theta(2s-p))+\frac{i}{2^p}\sum_{s=0}^{p}\binom{p}{s}\sin(\theta(2s-p))
=\frac{1}{2^{p-1}}\sum_{0\le{s}<p/2}\binom{p}{s}\cos((p-2s)\theta)+\frac{\delta_{p}}{2^{p}}\binom{p}{p/2}
\end{multline*}
 where $\delta_p=0$ for odd $p$, $\delta_p=1$  for even $p$, and the last equality follows by parity. Using this formula and employing the orthogonality of cosines (\ref{eq:cos-ortho}), we obtain for integer $p,k\ge0$:
\begin{multline*}
\int_0^{\pi}[\cos(\theta)]^p\cos(k\theta)d\theta
=\frac{1}{2^{p-1}}\sum_{0\le{s}<p/2}\binom{p}{s}\int_0^{\pi}\cos((p-2s)\theta)\cos(k\theta)d\theta
\\
+\frac{\delta_{p}}{2^{p}}\binom{p}{p/2}\int_0^{\pi}\cos(k\theta)d\theta
=\left\{\!\!\!\begin{array}{ll}0, &~\text{if}~p<k~\text{or}~p-k\equiv1\Mod2,\\\dfrac{\pi}{2^p}\displaystyle\binom{p}{(p-k)/2}, &~\text{if}~p-k\equiv0\Mod2.
\end{array}\right.
\end{multline*}
Substituting this formula into the definition of $A_l(k_1,k_2)$ we compute:
\begin{multline*}
A_l(k_1,k_2)=\int\limits_{0}^{\pi}\!\!\!\int\limits_{0}^{\pi}\sum_{s=0}^{l}\binom{l}{s}[\cos(\theta_1)]^s[\cos(\theta_1)]^{l-s}
\cos(k_1\theta_1)\cos(k_2\theta_2)d\theta_1d\theta_2
\\
=\sum_{s=0}^{l}\binom{l}{s}\int\limits_{0}^{\pi}[\cos(\theta_1)]^s\cos(k_1\theta_1)d\theta_1\int\limits_{0}^{\pi}[\cos(\theta_1)]^{l-s}\cos(k_2\theta_2)d\theta_2
\\
=\frac{\pi^2}{2^l}\sum_{\substack{k_1\le{s}\le{l-k_2}\\s-k_1=2m_1\\l-s-k_2=2m_2}}\binom{l}{s}\binom{s}{(s-k_1)/2}\binom{l-s}{(l-s-k_2)/2}
=\left\{\!\!\!\begin{array}{ll}0,~&\text{if}~l-k_1-k_2~\text{is odd}
\\
A_{2N+k_1+k_2}(k_1,k_2)~&\text{if}~l-k_1-k_2=2N.
\end{array}\right.
\end{multline*}
Further,
\begin{multline*}
A_{2N+k_1+k_2}(k_1,k_2)=\frac{\pi^2}{2^l}\sum_{\substack{0\le{n}\le{2N}\\n~\text{is even}}}\binom{2N+k_1+k_2}{n+k_1}\binom{n+k_1}{n/2}\binom{2N+k_2-n}{(2N-n)/2}
\\
=\frac{\pi^2}{2^l}\sum_{m=0}^{N}\binom{2N+k_1+k_2}{2m+k_1}\binom{2m+k_1}{m}\binom{2(N-m)+k_2}{N-m}=\frac{\pi^2}{2^l}S_{k_1,k_2}(N),
\end{multline*}
where the last equality is the definition of $S_{k_1,k_2}(N)$. It remains to calculate this binomial sum.  An application of Zeiberger's algorithm \cite[section 3.11]{AAR} (using, for instance, Fast Zeilberger Package by Peter Paule and Markus Schorn \cite{Paule-Schorn}) gives the recurrence:
$$
S_{k_1,k_2}(N+1)=\frac{(2N+k_1+k_2+1)^2(2N+k_1+k_2+2)^2}{(N+1)(N+k_1+1)(N+k_2+1)(N+k_1+k_2+1)}S_{k_1,k_2}(N).
$$
Using the initial value $S_{k_1,k_2}(0)=\binom{k_1+k_2}{k_1}$ we can solve this first order recurrence by straightforward back-substitution to obtain:
$$
S_{k_1,k_2}(N)=\binom{k_1+k_2}{k_1}\frac{(k_1+k_2+1)_{2N}^2}{(k_1+1)_{N}(k_2+1)_{N}(k_1+k_2+1)_{N}N!}.
$$
Plugging this formula into the expression for $A_{2N+k_1+k_2}(k_1,k_2)$ calculated above we arrive at (\ref{eq:A-computed}). $\hfill\square$
\par

\smallskip

Now we are in the position to give a proof of Theorem~\ref{th:Xi2representation}.

\medskip

\noindent\textbf{Proof of Theorem~\ref{th:Xi2representation}.} Using the generating function (\ref{eq:LegendreGenerating})
for Legendre's polynomials we have
$$
\frac{1}{\sqrt{1+x^2-x(\cos(\theta_1)+\cos(\theta_2))}}=\sum_{n=0}^{\infty}P_n((\cos(\theta_1)+\cos(\theta_2))/2)x^n
$$
The explicit formula for Legendre's polynomials \cite[3.6.1(15), 3.6.1(16)]{HTF1} reads
$$
P_n(y)=\frac{1}{2^n}\sum_{j=0}^{[n/2]}\frac{(-1)^j(2n-2j)!}{(n-2j)!(n-j)!j!}y^{n-2j}
=\left\{\!\!\!\begin{array}{l}
\dfrac{(-1)^m(2m)!}{4^m(m!)^2}{}_2F_1(-m,m+1/2;1/2;y^2),~~~~~n=2m,
\\[12pt]
\dfrac{(-1)^m(2m+1)!}{4^m(m!)^2}y{}_2F_1(-m,m+3/2;3/2;y^2),~~n=2m+1.
\end{array}
\right.
$$
Hence,
\begin{multline*}
4\Xi^{[\k]}_2(x)=\int_{0}^{\pi}\int_{0}^{\pi}
\frac{\cos(k_1\theta_1)\cos(k_2\theta_2)d\theta_1d\theta_2}{\sqrt{1+x^2-x(\cos(\theta_1)+\cos(\theta_2))}}
\\
=\sum_{m=0}^{\infty}x^{2m}\!\!\!\iint\limits_{0\le\theta_1,\theta_2\le\pi}\!\!\!\cos(k_1\theta_1)\cos(k_2\theta_2)P_{2m}((\cos(\theta_1)+\cos(\theta_2))/2)d\theta_1d\theta_2
\\
+\sum_{m=0}^{\infty}x^{2m+1}\!\!\!\iint\limits_{0\le\theta_1,\theta_2\le\pi}\!\!\!\cos(k_1\theta_1)\cos(k_2\theta_2)P_{2m+1}((\cos(\theta_1)+\cos(\theta_2))/2)d\theta_1d\theta_2
\\[7pt]
=\underbrace{\sum_{m=0}^{\infty}\frac{(-x^2)^{m}(2m)!}{4^m(m!)^2}
\sum_{j=0}^{m}\frac{(-m)_j(m+1/2)_j}{4^j(1/2)_jj!}
\!\!\!\iint\limits_{0\le\theta_1,\theta_2\le\pi}\!\!\!(\cos(\theta_1)+\cos(\theta_2))^{2j}\cos(k_1\theta_1)\cos(k_2\theta_2)d\theta_1d\theta_2}_{=T_{even}}
\\[10pt]
+\underbrace{x\sum_{m=0}^{\infty}\frac{(-x^2)^m(2m+1)!}{4^m(m!)^2}
\sum_{j=0}^{m}\frac{(-m)_j(m+3/2)_j}{2^{2j+1}(3/2)_jj!}
\!\!\!\iint\limits_{0\le\theta_1,\theta_2\le\pi}\!\!\!(\cos(\theta_1)+\cos(\theta_2))^{2j+1}\cos(k_1\theta_1)\cos(k_2\theta_2)
d\theta_1d\theta_2}_{T_{odd}}.
\end{multline*}

Assume first that $k_1+k_2$ is even and denote $k=(k_1+k_2)/2$.  Then by Lemma~\ref{lm:Al} we get $T_{odd}=0$. Further applying Lemma~\ref{lm:Al} to the even part we obtain:
\begin{multline*}
T_{even}=\sum_{m=0}^{\infty}\frac{(-x^2)^{m}(2m)!}{4^m(m!)^2}\sum_{j=0}^{m}\frac{(-m)_j(m+1/2)_j}{4^j(1/2)_jj!}A_{2j}(k_1,k_2)
\\
=\pi^2\binom{k_1+k_2}{k_1}\sum_{m=k}^{\infty}\frac{(-x^2)^{m}(2m)!}{4^m(m!)^2}\sum_{j=k}^{m}\frac{(-m)_j(m+1/2)_j(2k+1)_{2j-2k}^2}{4^{2j}(1/2)_j(k_1+1)_{j-k}(k_2+1)_{j-k}(k_1+k_2+1)_{j-k}(j-k)!j!}
\\
=\pi^2\binom{k_1+k_2}{k_1}\sum_{m=k}^{\infty}\frac{(-x^2)^{m}(1/2)_m}{m!}\sum_{i=0}^{m-k}\frac{(-m)_{i+k}(m+1/2)_{i+k}(2k+1)_{2i}^2}{4^{2(i+k)}(1/2)_{i+k}(k_1+1)_{i}(k_2+1)_{i}(2k+1)_{i}i!(i+k)!}
\\
=\frac{\pi^2(k+1/2)_kx^{2k}}{4^{2k}k!}\binom{k_1+k_2}{k_1}\sum_{n=0}^{\infty}\frac{(-x^2)^{n}(2k+1/2)_{n}}{n!}\sum_{i=0}^{n}\frac{(-n)_{i}(n+2k+1/2)_{i}(k+1/2)_i(k+1)_{i}}{(k_1+1)_{i}(k_2+1)_{i}(2k+1)_{i}i!}
\\
=\frac{\pi^2(k+1/2)_kx^{2k}}{4^{2k}k!}\binom{k_1+k_2}{k_1}\sum_{n=0}^{\infty}\frac{(-x^2)^{n}(2k+1/2)_{n}}{n!}
{_{4}F_3}\!\left(\!\!\begin{array}{l}-n,n+2k+1/2,k+1/2,k+1\\k_1+1,k_2+1,2k+1\end{array}\vline\,\,1\!\right)
\\
=\frac{\pi^2(1/2)_{2k}(2x)^{2k}}{4^{2k}(2k)!}\binom{k_1+k_2}{k_1}\sum_{n=0}^{\infty}\frac{(-x^2)^{n}(2k+1/2)_{n}}{n!}
{_{4}F_3}\!\left(\!\!\begin{array}{l}-n,n+2k+1/2,k+1/2,k+1\\k_1+1,k_2+1,2k+1\end{array}\vline\,\,1\!\right),
\end{multline*}
where we have used the formulas $(a)_{i+k}=(a)_k(a+k)_i$ and $(a)_{2i}=4^i(a/2)_i(a/2+1/2)_i$.  Similarly, if $k_1+k_2=2k+1$ is odd then $T_{even}=0$ and
\begin{multline*}
T_{odd}=x\sum_{m=0}^{\infty}\frac{(-x^2)^m(3/2)_m}{m!}\sum_{j=0}^{m}\frac{(-m)_j(m+3/2)_j}{2^{2j+1}(3/2)_jj!}A_{2j+1}(k_1,k_2)
\\
=x\pi^2\binom{k_1+k_2}{k_1}\sum_{m=k}^{\infty}\frac{(-x^2)^m(3/2)_m}{m}\sum_{j=k}^{m}\frac{(-m)_j(m+3/2)_j(k_1+k_2+1)_{2(j-k)}^2}{4^{2j+1}(3/2)_j(k_1+1)_{j-k}(k_2+1)_{j-k}(k_1+k_2+1)_{j-k}(j-k)!j!}
\\
=\frac{x\pi^2}{4^{2k+1}}\binom{k_1+k_2}{k_1}\sum_{n=0}^{\infty}\frac{(-x^2)^{n+k}(3/2)_{n+k}}{(n+k)!}\sum_{i=0}^{n}\frac{(-n-k)_{i+k}(n+k+3/2)_{i+k}(k_1+k_2+1)_{2i}^2}{4^{2i}(3/2)_{i+k}(k_1+1)_{i}(k_2+1)_{i}(k_1+k_2+1)_{i}i!(i+k)!}
\\
=\frac{\pi^2(k+3/2)_kx^{2k+1}}{4^{2k+1}}\binom{k_1+k_2}{k_1}\sum_{n=0}^{\infty}\frac{(-x^2)^{n}(3/2+2k)_n}{n!}\times
\\
\sum_{i=0}^{n}\frac{(-n)_{i}(n+2k+3/2)_{i}(k_1+k_2+1)_{2i}^2}{4^{2i}(3/2+k)_{i}(k_1+1)_{i}(k_2+1)_{i}(k_1+k_2+1)_{i}i!(i+k)!}
\\
=\frac{\pi^2(k+3/2)_kx^{2k+1}}{4^{2k+1}k!}\binom{k_1+k_2}{k_1}\sum_{n=0}^{\infty}\frac{(-x^2)^{n}(2k+3/2)_n}{n!}
\sum_{i=0}^{n}\frac{(-n)_{i}(n+2k+3/2)_{i}(k+1)_i(k+3/2)_i}{(k_1+1)_{i}(k_2+1)_{i}(2k+2)_{i}i!}
\\
=\frac{\pi^2(k+3/2)_kx^{2k+1}}{4^{2k+1}k!}\binom{k_1+k_2}{k_1}\sum_{n=0}^{\infty}\frac{(-x^2)^{n}(2k+3/2)_n}{n!}
{_{4}F_3}\!\left(\!\!\begin{array}{l}-n,n+2k+3/2,k+3/2,k+1\\k_1+1,k_2+1,2k+2\end{array}\vline\,\,1\!\right)
\\
=\frac{\pi^2(1/2)_{2k+1}(2x)^{2k+1}}{4^{2k+1}(2k+1)!}\binom{k_1+k_2}{k_1}\sum_{n=0}^{\infty}\frac{(-x^2)^{n}(2k+3/2)_n}{n!}
{_{4}F_3}\!\left(\!\!\begin{array}{l}-n,n+2k+3/2,k+3/2,k+1\\k_1+1,k_2+1,2k+2\end{array}\vline\,\,1\!\right).
\end{multline*}
Substituting back the formulas for $k$ in each case we see that the expressions for  $T_{even}$ and $T_{odd}$ coincide!  Hence, in all cases
\begin{multline*}
4\Xi^{[\k]}_2(x)
=\frac{\pi^2(1/2)_{k_1+k_2}(2x)^{k_1+k_2}}{4^{k_1+k_2}k_1!k_2!}\sum_{n=0}^{\infty}\frac{(-x^2)^{n}(k_1+k_2+1/2)_n}{n!}
\\
\times
{_{4}F_3}\!\left(\!\!\begin{array}{l}-n,n+k_1+k_2+1/2,(k_1+k_2+2)/2,(k_1+k_2+1)/2\\k_1+1,k_2+1,k_1+k_2+1\end{array}\vline\,\,1\!\right).
\end{multline*}
This expression can be further reduced to a finite sum of ${}_3F_2$ hypergeometric functions.  To this end note that ${}_4F_3(1)$ appearing in the above formula is balanced (or Saalsch\"{u}tzian), which means that the sum of top parameters is one less than the sum of the bottom parameters.
For balanced ${}_4F_3(1)$ the next transformation is well-known \cite[Theorem~3.3.3]{AAR}:
$$
{_{4}F_3}\!\left(\!\!\begin{array}{l}-n,A,B,C\\E,F,G\end{array}\vline\,\,1\!\right)=\frac{(F-C)_n(G-C)_n}{(F)_n(G)_n}{_{4}F_3}\!\left(\!\!\begin{array}{l}-n,E-A,E-B,C\\E,E+F-A-B,E+G-A-B\end{array}\vline\,\,1\!\right)
$$
(provided that $-n+A+B+C+1=E+F+G$).  Assuming, without loss of generality, that $k_1-k_2=-l\le0$ and applying this formula with $A=(k_1+k_2+1)/2$, $B=(k_1+k_2+2)/2$, $C=n+k_1+k_2+1/2$, $E=k_1+1$, $F=k_2+1$, $G=k_1+k_2+1$ we get
\begin{multline*}
{_{4}F_3}\!\left(\!\!\begin{array}{l}-n,n+k_1+k_2+1/2,(k_1+k_2+1)/2,(k_1+k_2+2)/2\\k_1+1,k_2+1,k_1+k_2+1\end{array}\vline\,\,1\!\right)
\\
=\frac{(1/2)_n(k_1+1/2)_n}{(k_2+1)_n(k_1+k_2+1)_n}
{_{4}F_3}\!\left(\!\!\begin{array}{l}-n,n+k_1+k_2+1/2,(k_1-k_2+1)/2,(k_1-k_2)/2\\1/2,k_1+1/2,k_1+1\end{array}\vline\,\,1\!\right),
\end{multline*}
where we applied $(1/2-n+\alpha)_n=(-1)^n(1/2+\alpha)_n$ with $\alpha=0$ and $\alpha=k_1$.
Hence, for $l=2r$ we get by exchanging the order of summations and using $(-n)_j=(-1)^jn!/(n-j)!$ and $(a)_{n+j}=(a)_j(a+n)_j$:
\begin{multline*}
T_{even}=\frac{\pi^2(1/2)_{k_1+k_2}(2x)^{k_1+k_2}}{4^{k_1+k_2}k_1!k_2!}\sum_{n=0}^{\infty}\frac{(-x^2)^{n}(k_1+k_2+1/2)_{n}(1/2)_n(k_1+1/2)_n}{(k_2+1)_n(k_1+k_2+1)_nn!}
\\
\times
{_{4}F_3}\!\left(\!\!\begin{array}{l}-r,1/2-r,-n,n+k_1+k_2+1/2\\1/2,k_1+1/2,k_1+1\end{array}\vline\,\,1\!\right)
\\
=\frac{\pi^2(1/2)_{k_1+k_2}(2x)^{k_1+k_2}}{4^{k_1+k_2}k_1!k_2!}\sum_{n=0}^{\infty}\frac{(-x^2)^{n}(k_1+k_2+1/2)_{n}(1/2)_n(k_1+1/2)_n}{(k_2+1)_n(k_1+k_2+1)_nn!}
\\
\times\sum_{j=0}^{r}\frac{(-r)_j(1/2-r)_j(-n)_j(n+k_1+k_2+1/2)_j}{(1/2)_j(k_1+1/2)_j(k_1+1)_jj!}
\\
=\frac{\pi^2(1/2)_{k_1+k_2}(2x)^{k_1+k_2}}{4^{k_1+k_2}k_1!k_2!}\sum_{j=0}^{r}\frac{(-1)^j(-r)_j(1/2-r)_j(k_1+k_2+1/2)_{j}}{(1/2)_j(k_1+1/2)_j(k_1+1)_jj!}
\\
\times\sum_{n=j}^{\infty}\frac{(-x^2)^{n}(k_1+k_2+1/2+j)_n(1/2)_n(k_1+1/2)_n}{(k_2+1)_n(k_1+k_2+1)_n(n-j)!}
\\
=\frac{\pi^2(1/2)_{k_1+k_2}(2x)^{k_1+k_2}}{4^{k_1+k_2}k_1!k_2!}\sum_{j=0}^{r}\frac{(-r)_j(1/2-r)_j(k_1+k_2+1/2)_{2j}x^{2j}}
{(k_1+1)_j(k_2+1)_{j}(k_1+k_2+1)_{j}j!}
\\
\times\sum_{m=0}^{\infty}\frac{(-x^2)^{m}(k_1+k_2+1/2+2j)_{m}(1/2+j)_{m}(k_1+1/2+j)_{m}}{(k_2+1+j)_{m}(k_1+k_2+1+j)_{m}m!}
\\
=\frac{\pi^2(1/2)_{k_1+k_2}(2x)^{k_1+k_2}}{4^{k_1+k_2}k_1!k_2!}\sum_{j=0}^{r}\frac{(-r)_j(1/2-r)_j(k_1+k_2+1/2)_{2j}x^{2j}}{(k_1+1)_j(k_2+1)_{j}(k_1+k_2+1)_{j}j!}
\\
\times{_{3}F_2}\!\left(\!\!\begin{array}{l}1/2+j,k_1+1/2+j,k_1+k_2+1/2+2j\\k_2+1+j,k_1+k_2+1+j\end{array}\vline\,\,-x^2\!\right).
\end{multline*}
Similarly, if $k_2-k_1=l=2r+1$, $r=0,1,\ldots$ for $T_{odd}$ we get:
\begin{multline*}
T_{odd}=\frac{\pi^2(1/2)_{k_1+k_2}(2x)^{k_1+k_2}}{4^{k_1+k_2}k_1!k_2!}\sum_{n=0}^{\infty}\frac{(-x^2)^{n}(k_1+k_2+1/2)_{n}(1/2)_n(k_1+1/2)_n}{(k_2+1)_n(k_1+k_2+1)_nn!}
\\
\times
\frac{(1/2)_n(k_1+1/2)_n}{(k_2+1)_n(k_1+k_2+1)_n}
{_{4}F_3}\!\left(\!\!\begin{array}{l}-r,-1/2-r,-n,n+k_1+k_2+1/2\\1/2,k_1+1/2,k_1+1\end{array}\vline\,\,1\!\right)
\\
=\frac{\pi^2(1/2)_{k_1+k_2}(2x)^{k_1+k_2}}{4^{k_1+k_2}k_1!k_2!}\sum_{j=0}^{r}\frac{(-r)_j(-1/2-r)_j(k_1+k_2+1/2)_{2j}x^{2j}}{(k_1+1)_j(k_2+1)_{j}(k_1+k_2+1)_{j}j!}
\\
\times{_{3}F_2}\!\left(\!\!\begin{array}{l}1/2+j,k_1+1/2+j,k_1+k_2+1/2+2j\\k_2+1+j,k_1+k_2+1+j\end{array}\vline\,\,-x^2\!\right).
\end{multline*}
Uniting both formulas we finally arrive at (\ref{eq:theta2-3F2}).$\hfill\square$

\begin{corollary}\label{cr:dif-dif}
The super- and sub- diagonal functions $\Xi_2^{[k,k+1]}=\Xi_2^{[k+1,k]}$ satisfy the following differential-difference relation:
$$
\frac{1}{2}\left(\frac{1}{x}+(1-x^2)\frac{d}{dx}\right)\Xi_2^{[k,k+1]}(x)=(k+1/4)\Xi_2^{[k,k]}(x)-(k+3/4)\Xi_2^{[k+1,k+1]}(x).
$$
\end{corollary}
The proof will require the following lemma.
\begin{lemma}\label{lm:3F2-contiguous}
The following contiguous relation holds for ${}_3F_{2}$ \emph{(}here $\delta=y\dfrac{d}{dy}$\emph{):}
\begin{multline}\label{eq:3F2-contiguous}
[b_1b_2+a_2(a_3-a_1)y+(b_2+(a_3-a_1)y)\delta]
{_{3}F_2}\!\left(\!\!\begin{array}{l}a_1,a_2,a_3\\b_1+1,b_2+1\end{array}\vline\,\,y\!\right)
\\
=\frac{a_2a_3(b_2-a_1+1)y}{b_2+1}{_{3}F_2}\!\left(\!\!\begin{array}{l}a_1,a_2+1,a_3+1\\b_1+1,b_2+2\end{array}\vline\,\,y\!\right)
+(b_1b_2){_{3}F_2}\!\left(\!\!\begin{array}{l}a_1,a_2,a_3-1\\b_1,b_2\end{array}\vline\,\,y\!\right).
\end{multline}
\end{lemma}
\textbf{Proof.}  The proof hinges on the next two identities:
$$
\left[1-\frac{a_1a_2}{b_1b_2}y+\frac{b_1+b_2-(a_1+a_2)y}{b_1b_2}\delta+\frac{1-y}{b_1b_2}\delta^2\right]
\!{_{3}F_2}\!\left(\!\!\begin{array}{l}a_1,a_2,a_3\\b_1+1,b_2+1\end{array}\!\vline\,y\!\right)
\!=\!{_{3}F_2}\!\left(\!\!\begin{array}{l}a_1,a_2,a_3-1\\b_1,b_2\end{array}\!\vline\,y\!\right)
$$
and
$$
\left[1-\frac{b_1-(a_2+a_3)y}{a_2a_3}\delta
-\frac{(1-y)}{a_2a_3}\delta^2\right]
\!{_{3}F_2}\!\left(\!\!\begin{array}{l}a_1,a_2,a_3\\b_1+1,b_2+1\end{array}\!\vline\,y\!\right)
\!=\!\frac{b_2-a_1+1}{b_2+1}\!{_{3}F_2}\!\left(\!\!\begin{array}{l}a_1,a_2+1,a_3+1\\b_1+1,b_2+2\end{array}\!\vline\,y\!\right).
$$
These identities can be obtained from formulas (15), (26) found in \cite[Table~5]{Sasaki} by simple renaming of variables or by using the Mathematica package HYPERDIRE \cite{BKK}.
The claim now follows by removing the term containing
$$
\delta^2\!{_{3}F_2}\!\left(\!\!\begin{array}{l}a_1,a_2,a_3\\b_1+1,b_2+1\end{array}\!\vline\,y\!\right)
$$
from both equations. $\hfill\square$

\smallskip

\noindent\textbf{Proof of Corollary~\ref{cr:dif-dif}.}
Using (\ref{eq:corollaryXi2}) define
\begin{equation}\label{eq:hattheta}
\hat{\Xi}_{k_1,k_2}(x)=\frac{4^{k_1+k_2+1}k_1!k_2!}{\pi^2(1/2)_{k_1+k_2}(2x)^{k_1+k_2}}\Xi_2^{[k_1,k_2]}(x)
={_{3}F_2}\!\left(\!\!\begin{array}{l}1/2,k_1+1/2,k_1+k_2+1/2\\k_2+1,k_1+k_2+1\end{array}\vline\,\,-x^2\!\right).
\end{equation}
Denoting $a_1=1/2$, $a_2=k_1+1/2$, $a_3=k_1+k_2+1/2$, $b_1=k_2$, $b_2=k_1+k_2$, we have
\begin{equation}\label{eq:Dhattheta}
\frac{x}{2}\frac{d}{dx}\hat{\Xi}_{k_1,k_2}(x)=\frac{a_1a_2a_3(-x^2)}{(b_1+1)(b_2+1)}{_{3}F_2}\!\left(\!\!\begin{array}{l}a_1+1,a_2+1,a_3+1\\b_1+2,b_2+2\end{array}\vline\,\,-x^2\!\right).
\end{equation}
Next, from (\ref{eq:3F2-contiguous}) by an application of the differential operator $\delta$ we have
\begin{multline*}
(b_1b_2+a_2(a_3-a_1)y){_{3}F_2}\!\left(\!\!\begin{array}{l}a_1,a_2,a_3\\b_1+1,b_2+1\end{array}\vline\,\,y\!\right)
\\
+(b_2+(a_3-a_1)y)\frac{a_1a_2a_3y}{(b_1+1)(b_2+1)}{_{3}F_2}\!\left(\!\!\begin{array}{l}a_1+1,a_2+1,a_3+1\\b_1+2,b_2+2\end{array}\vline\,\,y\!\right)
\\
=\frac{a_2a_3(b_2-a_1+1)y}{b_2+1}{_{3}F_2}\!\left(\!\!\begin{array}{l}a_1,a_2+1,a_3+1\\b_1+1,b_2+2\end{array}\vline\,\,y\!\right)
+(b_1b_2){_{3}F_2}\!\left(\!\!\begin{array}{l}a_1,a_2,a_3-1\\b_1,b_2\end{array}\vline\,\,y\!\right).
\end{multline*}
Writing $y=-x^2$ and employing (\ref{eq:hattheta}) and  (\ref{eq:Dhattheta}) we get
\begin{multline*}
(b_1b_2+a_2(a_3-a_1)(-x^2))\hat{\Xi}_{k_1,k_2}(x)
+(b_2+(a_3-a_1)(-x^2))\frac{x}{2}\frac{d}{dx}\hat{\Xi}_{k_1,k_2}(x)
\\
=\frac{a_2a_3(b_2-a_1+1)(-x^2)}{b_2+1}\hat{\Xi}_{k_1+1,k_2}(x)
+(b_1b_2)\hat{\Xi}_{k_1,k_2-1}(x).
\end{multline*}
or, substituting back the definition of parameters $a_1$, $a_2$, $a_3$, $b_1$, $b_2$, we obtain:
\begin{multline*}
(k_2-x^2(k_1+1/2))\hat{\Xi}_{k_1,k_2}(x)
+(1-x^2)\frac{x}{2}\frac{d}{dx}\hat{\Xi}_{k_1,k_2}(x)
\\
=k_2\hat{\Xi}_{k_1,k_2-1}(x)-x^2\frac{(k_1+1/2)(k_1+k_2+1/2)^2}{(k_1+k_2)(k_1+k_2+1)}\hat{\Xi}_{k_1+1,k_2}(x).
\end{multline*}
Recalling that $k_2=k_1+1=k+1$ this amounts to
\begin{multline*}
(k+1-x^2(k+1/2))\hat{\Xi}_{k,k+1}(x)
+(1-x^2)\frac{x}{2}\frac{d}{dx}\hat{\Xi}_{k,k+1}(x)
\\
=(k+1)\hat{\Xi}_{k,k}(x)-x^2\frac{(k+3/4)^2}{(k+1)}\hat{\Xi}_{k+1,k+1}(x).
\end{multline*}
Finally, substituting the definition of $\hat{\Xi}_{k_1,k_2}$ form (\ref{eq:hattheta}) we get:
\begin{multline*}
2(k+1-x^2(k+1/2))\Xi_2^{[k,k+1]}(x)+(1-x^2)x^{2k+2}\frac{d}{dx}\frac{1}{x^{2k+1}}\Xi_2^{[k,k+1]}(x)
\\
=(1/2+2k)x\Xi_2^{[k,k]}(x)-\frac{4x(k+3/4)^2}{(2k+3/2)}\Xi_2^{[k+1,k+1]}(x).
\end{multline*}
Taking derivative we finally arrive at:
$$
\left[1+x(1-x^2)\frac{d}{dx}\right]\Xi_2^{[k,k+1]}(x)=2x[(k+1/4)\Xi_2^{[k,k]}(x)-(k+3/4)\Xi_2^{[k+1,k+1]}(x)].~~~~~\square
$$
\par
We can now present a result similar to Corollary~\ref{cr:rank1diffeq} in case $N=2$.
\begin{corollary}\label{cr:rank2diffeq}
Suppose in addition to  the assumptions of Theorem~\ref{th:Xi2representation} that $k_2-k_1\le 1$. Then
the function $\Xi^{[{\bf k}]}_2(x)$ obeys the third-order differential equation
\begin{multline}\label{eq:Xi2diffeq}
\Bigl\{x^3(1+x^2)D_x^3+x^2((k_1-k_2+6)(x^2-1)+9)D_x^2-x((k_1+k_2)^2(x^2+1)+(k_2-k_1-7/3)(3x^2-1)-10/3)D_x
\\
+x^2((k_1+k_2)^2-1)(k_2-k_1-1)+(k_1-k_2)(k_1+k_2)^2\Bigr\}\Xi^{{\bf [k]}}_2(x)=0,
\end{multline}
where $D_x=d/dx$.
\end{corollary}
\textbf{Proof.} Begin with the differential equation \cite[4.2(2)]{HTF1}
\begin{equation}\label{eq:3F2diffeq}
\Bigl[{\cal T}({\cal T}+b_1-1)({\cal T}+b_2-1)-z({\cal T}+a_1)({\cal T}+a_2)({\cal T}+a_3)\Bigr]f(z)=0
\end{equation}
satisfied by the function
$$
f(z)={_{3}F_2}\!\left(\!\!\begin{array}{l}a_1,a_2,a_3\\b_1,b_2\end{array}\vline\,\,z\!\right),
$$
where the operator ${\cal T}$ is defined by
$$
{\cal T}:=z\frac{d}{dz}.
$$
Next to find the equation satisfied by $f(-x^2)$ we apply the easily verifiable identities:
$$
[{\cal T}f(z)]_{\vert z=-x^2}=\frac{x}{2}D_xf(-x^2), ~~~~[{\cal T}^2f(z)]_{\vert z=-x^2}=\frac{1}{4}\left[x^2D_x^2+xD_x\right]f(-x^2),
$$
$$
[{\cal T}^3f(z)]_{\vert z=-x^2}=\frac{1}{8}\left[x^3D_x^3+3x^2D_x^2+xD_x\right]f(-x^2).
$$
Substituting these expressions into (\ref{eq:3F2diffeq}) we obtain:
\begin{multline*}
\Bigl\{\frac{x^3}{8}(1+x^2)D_x^3+\frac{x^2}{4}(x^2(a_1+a_2+a_3+3/2)+b_1+b_2-1/2)D_x^2
\\
+\frac{x^2}{4}(x^2(a_1+a_2+a_3+2a_1a_2+2a_1a_3+2a_2a_3+1/2)+(2b_1-1)(2b_2-1)/2)D_x+x^2a_1a_2a_3\Bigr\}f(-x^2)=0.
\end{multline*}
Finally, to find the differential equation satisfied by $x^{\alpha}f(-x^2)$ we utilize the straightforward formulas
\begin{align*}
&D_x[g(x)]=x^{-\alpha}\left\{D_x-\alpha{x^{-1}}\right\}x^{\alpha}g(x),
\\[5pt]
&D_x^2[g(x)]=x^{-\alpha}\left\{D_x^2-2\alpha{x^{-1}}D_x+\alpha(\alpha+1)x^{-2}\right\}x^{\alpha}g(x),
\\[5pt]
&D_x^3[g(x)]=x^{-\alpha}\left\{D_x^3-3\alpha{x^{-1}}D_x^2+3\alpha(\alpha+1)x^{-2}D_x-\alpha(\alpha^2+3\alpha+2)x^{-3} \right\}x^{\alpha}g(x).
\end{align*}
Substituting these formulas in the above differential equation with $\alpha=k_1+k_2$, $g(x)=f(-x^2)$, $a_1=1/2$, $a_2=k_1+1/2$, $a_2=k_1+k_2+1/2$, $b_1=k_2+1$, $b_2=k_1+k_2+1$
after substantial amount of calculations and simplifications we arrive at (\ref{eq:Xi2diffeq}). $\hfill\square$

Let us note, that for the ''lowest mode'' $\Xi$-function, i.e. in case $k_1=k_2=0$, Corollary~\ref{cr:rank2diffeq} implies
$$
\Bigl\{x^3(1+x^2)D_x^3+3x^2(2x^2+1)D_x^2+x(7x^2+1)D_x+x^2\Bigr\}\Xi^{[0,0]}_2(x)=0.
$$

\bigskip

\section{Conclusions and further tasks}

We have expressed functions $\Xi_1^{[k]}$ and $\Xi_2^{[{\bf k}]}$ in terms of the Gauss and Clausen hypergeometric functions $_{2}F_1$ and $_{3}F_2$, respectively. We have also found
differential equations and differential-difference relations satisfied by $\Xi$-functions. We believe that the higher rank ($N\ge3)$ $\Xi_N^{[{\bf k}]}(x)$ functions could be expressed in terms of (probably multivariate) hypergeometric functions.

\par
A study of the higher rank $\Xi$-functions is more of a mathematical interest and will be performed at a later stage. Further, we plan to consider special functions appearing from similar problems of the analysis of the operators $K_\varphi$ in other physically important three-dimensional domains, e.g. cylindrical and spherical. For example, such functions for a constant field $\varphi$ in cylindrical domains \cite{MelnikovJMP2017-2} are
$$
\psi^{[k]}_\zeta(x)=\int_0^{2\pi}\!\!\!\frac{\cos(k\theta)d\theta}{\sqrt{1+\zeta^2+x^2-2x\cos\theta}}.
$$

\end{document}